\documentclass[11pt]{article}

\usepackage{mathpazo}   
\usepackage{tikz}

\usepackage{nameref}
\usepackage{empheq}
\usepackage{comment}
\usepackage[shortlabels,inline]{enumitem}
\usepackage[colorlinks=true,
linkcolor=refkey,
urlcolor=lblue,
citecolor=red]{hyperref}
\usepackage[doc,wmm,hhb]{optional}
\usepackage{xcolor}

\usepackage{float}
\usepackage{soul}
\usepackage{graphicx}
\definecolor{labelkey}{rgb}{0,0.08,0.45}
\definecolor{refkey}{rgb}{0,0.6,0.0}
\definecolor{Brown}{rgb}{0.45,0.0,0.05}
\definecolor{lime}{rgb}{0.00,0.8,0.0}
\definecolor{lblue}{rgb}{0.5,0.5,0.99}
\definecolor{OliveGreen}{rgb}{0,0.6,0}
\definecolor{tyrianpurple}{rgb}{0.4, 0.01, 0.24}

\colorlet{hlcyan}{cyan!30}

\usepackage{stmaryrd}
\usepackage{amssymb}

\makeatletter
\def\namedlabel#1#2{\begingroup
   \def\@currentlabel{#2}%
   \label{#1}\endgroup
}
\makeatother

\newcommand{\seppthree}{\setlength{\itemsep}{-3pt}}

\usepackage[margin=0.92in,footskip=0.25in]{geometry}
\newcommand{\nnn}{\ensuremath{{n\in{\mathbb N}}}}

\newcommand{\menge}[2]{\big\{{#1}~\big |~{#2}\big\}}

\newcommand{\fenv}[1]%
{\ensuremath{\,\overrightarrow{\operatorname{env}}_{#1}}}
\newcommand{\benv}[1]%
{\ensuremath{\,\overleftarrow{\operatorname{env}}_{#1}}}

\newcommand{\scal}[2]{\left\langle{#1},{#2}  \right\rangle}

\newcommand{\RR}{\ensuremath{\mathbb R}}
\newcommand{\ran}{\ensuremath{{\operatorname{ran}}\,}}

\newcommand{\conv}{\ensuremath{\operatorname{conv}\,}}

\newcommand{\Fix}{\ensuremath{\operatorname{Fix}}}

\newcommand{\sperp}{{\scriptscriptstyle\perp}}

{\begin{list}{}{%
\settowidth{\labelwidth}{\textrm{#1~}}%
\setlength{\leftmargin}{\labelwidth+\labelsep}}}
{\end{list}}
\usepackage{amsthm}
\makeatletter
\def\th@plain{%
	\thm@notefont{}
	\itshape 
}
\def\th@definition{%
	\thm@notefont{}
	\normalfont 
}
\makeatother
\usepackage[capitalize,nameinlink]{cleveref}
\crefname{equation}{}{equations}
\crefname{chapter}{Appendix}{chapters}
\crefname{item}{}{items}
\crefname{enumi}{}{}
\newtheorem{theorem}{Theorem}[section]

\newtheorem{corollary}[theorem]{Corollary}

\newtheorem{proposition}[theorem]{Proposition}

\newtheorem{example}[theorem]{Example}
\newtheorem{ex}[theorem]{Example}

\newtheorem{remark}[theorem]{Remark}

\providecommand{\ds}{\displaystyle}

\providecommand{\norm}[1]{\lVert#1\rVert}

\providecommand{\RR}{\mathbb{R}}

\providecommand{\conv}{\operatorname{conv}}
\providecommand{\ran}{\operatorname{ran}}

\providecommand{\ran}{\operatorname{ran}}

\providecommand{\RR}{\mathbb{R}}

\definecolor{myblue}{rgb}{0.9,0.9,0.98}

  \newcommand*\mybluebox[1]{%
    \colorbox{myblue}{\hspace{1em}#1\hspace{1em}}}

\allowdisplaybreaks 

\begin{document}

\setlength{\abovedisplayskip}{8pt}
\setlength{\belowdisplayskip}{8pt}	
	

%

\author{
Heinz H.\ Bauschke\thanks{
Department of Mathematics, University
of British Columbia,
Kelowna, B.C.\ V1V~1V7, Canada. E-mail:
\texttt{heinz.bauschke@ubc.ca}.},~~~
Dayou Mao\thanks{
Department of Combinatorics and Optimization, 
University of Waterloo,
Waterloo, Ontario N2L~3G1, Canada.
  E-mail: \texttt{d6mao@uwaterloo.ca}.}
~~~and~
Walaa M.\ Moursi\thanks{
Department of Combinatorics and Optimization, 
University of Waterloo,
Waterloo, Ontario N2L~3G1, Canada.
  E-mail: \texttt{walaa.moursi@uwaterloo.ca}.}
}

\title{\textsf{
How to project onto the intersection of\\ a closed affine subspace and a hyperplane 
}
}

\date{June 22, 2022}

\maketitle

\begin{abstract}
Let $A$ be a closed affine subspace and let $B$ be a hyperplane in a Hilbert space. Suppose we are given their associated nearest point mappings $P_A$ and $P_B$, respectively. We present a formula for the projection onto their intersection $A\cap B$. 
As a special case, we derive a formula for the projection onto the intersection of two 
hyperplanes. 
These formulas provides useful information even if $A\cap B$ is empty.  Examples and numerical experiments are also provided.
\end{abstract}

{ 
\small
\noindent
{\bfseries 2020 Mathematics Subject Classification:}
{Primary 
15A04, 
41A50, 
47A50;
Secondary 
47H09,
90C25. 
}

\noindent {\bfseries Keywords:}
affine subspace, 
Hilbert space,
hyperplane,
linear equation, 
projection onto an intersection. 
}

\section{Introduction}

Throughout, we assume that 
\begin{empheq}[box=\mybluebox]{equation}
\text{$X$ is
a real Hilbert space with inner product 
$\scal{\cdot}{\cdot}\colon X\times X\to\RR$, }
\end{empheq}
and induced norm $\|\cdot\|$.
We also assume that 
\begin{empheq}[box=\mybluebox]{equation}
\label{e:A}
\text{$A$ is a closed affine subspace of $X$ with 
parallel space $U := A-A$.}
\end{empheq}
Denote the nearest point mapping associated with $A$ by $P_A$, and 
set $a_0 := P_A(0)$. Then $a_0\in U^\perp$ and
$P_A(x)=a_0 + P_U(x)$. This formula allows us to move back and forth 
between $P_A$ and $P_U$ as needed. 
Furthermore, we assume that $B$ is a hyperplane given by 
\begin{empheq}[box=\mybluebox]{equation}
\label{e:B}
\text{$B := \menge{x\in X}{\scal{x}{v}=\beta}$, where $v\in X$ and $\|v\|=1$, }
\end{empheq}
which in turn yields 
$B = \beta v+\{v\}^\perp$ and $P_B\colon x\mapsto x-(\scal{x}{v}-\beta)v$. 

\emph{The aim of this note is to present and prove a formula for $P_{A\cap B}$ that 
relies on $P_A$ and $P_B$.} 

Indeed, we have: 

\begin{theorem}[main result]
\label{t:main}
Suppose that $A\cap B\neq\varnothing$, and let $x\in X$.
Then exactly one of the following holds:
\begin{enumerate}
\item \label{t:main1}
If $P_U(v)=0$, then $P_{A\cap B}(x)=P_A(x)$.
\item \label{t:main2}
If $P_U(v)\neq 0$, then 
$\displaystyle P_{A\cap B}(x) = P_A(x) + \frac{\beta-\scal{P_A(x)}{v}}{\|P_U(v)\|^2}P_U(v)$. 
\end{enumerate}
\end{theorem}

The analysis we carry out will reveal useful information \emph{even if $A\cap B=\varnothing$.}

\emph{Organization of the paper.}
In \cref{sec:zero}, we analyze the case when $P_U(v)=0$. 
In \cref{sec:nonzero}, we turn to the situation when $P_U(v)\neq 0$.
A utility version of our main result as well as the special case of
two hyperplanes is discussed in \cref{sec:uti}.
In \cref{sec:numexp}, we outline a numerical application. 

\section{The case when $P_U(v)=0$}
\label{sec:zero}

Throughout this section, we assume that 
\begin{empheq}[box=\mybluebox]{equation}
P_U(v)=0.
\end{empheq}
Note that this assumption allows the following additional characterizatons: 
\begin{equation}
\label{e:addchar}
\text{$v\in U^\perp$
$\Leftrightarrow$
$\RR v\subseteq U^\perp$
$\Leftrightarrow$
$U \subseteq \{v\}^\perp$.}
\end{equation}
We shall not impose that $A\cap B\neq\varnothing$.
In order to handle this case, we well define 
\begin{empheq}[box=\mybluebox]{equation}
g := P_{{B-A}}(0)
\;\;\text{and}\;\;
E := A\cap (B-g)\neq\varnothing.
\end{empheq}
Note that $B-A = (U+\{v\}^\perp) + (\beta v-a_0)$.
Because $\{v\}^\perp$ is codimension $1$, the Minkowski sum
$U+\{v\}^\perp$ is closed by combining \cite[Theorem~9.35 and Corollary~9.37]{Deutsch}. Hence $B-A$ is a closed affine subspace which makes the vector $g$ not only well-defined but it also yields $E\neq\varnothing$. Moreover, 
\cite[Example~2.2]{01} yields
\begin{equation}
\label{e:g}
g = P_{U^\perp\cap\RR v}(b-a)\in U^\perp \cap \RR v
\end{equation}
no matter how $(a,b)\in A\times B$ is chosen. 
Note that 
\begin{equation}
\text{if $A\cap B\neq\varnothing$, then $g=0$ and so $E = A\cap B$;}
\end{equation}
consequently, $E$ can be thought of as the 
\emph{generalized} intersection of $A$ and $B$. 
Finally, by \cite[Lemma~2.2.(i)]{02}, 
the generalized intersection $E$ has also a description that does not involve  the gap vector $g$: 
\begin{equation}
\label{e:altE}
E = \Fix(P_AP_B) := \menge{x\in X}{x=P_A(P_B(x))}. 
\end{equation}

We now prove

\begin{theorem}
\label{p:1}
The assumption that $P_U(v)=0$ implies 
$A \subseteq B-g$ and so $E=A$.
\end{theorem}
\begin{proof}
Take $e\in E = A\cap (B-g)$.
Then $e=P_A(e)=a_0+P_U(e)$ and 
$e=\beta v+v^\sperp-g$ for some $v^\sperp\in \{v\}^\perp$. 
Using the linearity of $P_{U^\perp}$, the fact that $a_0\in U^\perp$, 
\cref{e:addchar}, and \cref{e:g}, we obtain 
$a_0 = P_{U^\perp}(e)=
\beta P_{U^\perp}(v)+P_{U^\perp}(v^\sperp) - P_{U^\perp}(g)
=\beta v + v^\sperp - P_{U}(v^\sperp) - g$. 
Therefore, using \cref{e:addchar} again,
\begin{subequations}
\begin{align}
A &= a_0 + U = \beta v + v^\sperp - P_{U}(v^\sperp) - g + U\\
&=\beta v + v^\sperp +(U-P_U(v^\sperp)) - g\\
&=\beta v + v^\sperp +U - g\\
&\subseteq \beta v + \{v\}^\perp + \{v\}^\perp -g\\
&= (\beta v + \{v\}^\perp) -g\\
&= B-g,
\end{align}
\end{subequations}
as claimed. 
Because $E=A\cap (B-g)$, it now follows that $E=A$. 
\end{proof}

\begin{corollary}
\label{c:1}
The assumption that $P_U(v)=0$ yields $P_E=P_A$. 
Let $x\in X$. 
Then exactly one of the following holds:
\begin{enumerate}
\item \label{c:1.1}
$P_A(x)\in B$, $g=0$, $E=A\cap B\neq \varnothing$, 
and $P_{A\cap B}(x) = P_A(x)$.
\item \label{c:1.2}
$P_A(x)\notin B$, $g\neq 0$, $E\neq A\cap B=\varnothing$, and 
$P_E(x)=P_A(x)$.
\end{enumerate}
\end{corollary}
\begin{proof}
This is a direct consequence of \cref{p:1}. 
\end{proof}

\begin{remark}
Note that \cref{c:1}\cref{c:1.1} yields \cref{t:main}\cref{t:main1}. 
\end{remark}

\section{The case when $P_U(v)\neq 0$}

\label{sec:nonzero}

Throughout this section, we assume that 
\begin{empheq}[box=\mybluebox]{equation}
P_U(v)\neq 0.
\end{empheq}
Then
\begin{equation}
\label{e:220616b}
0< \|P_U(v)\|^2 = \scal{P_U(v)}{v}. 
\end{equation}
Now set
\begin{equation}
\label{e:Q}
Q\colon X\to X\colon 
x\mapsto P_A(x) + \frac{\beta-\scal{P_A(x)}{v}}{\|P_U(v)\|^2}P_U(v) 
\in A+ U = A.
\end{equation}

\begin{proposition}
\label{p:neq1}
The assumption that $P_U(v)\neq 0$ implies 
$\ran Q\subseteq A\cap B$; 
hence, $A\cap B\neq\varnothing$. 
\end{proposition}
\begin{proof}
Let $x\in X$. 
Using \cref{e:Q} and \cref{e:220616b}, we have 
$Q(x)\in A$ and
\begin{subequations}
\begin{align}
\scal{Q(x)}{v}
&= \scal{P_A(x)}{v}
+ \frac{\beta-\scal{P_A(x)}{v}}{\|P_U(v)\|^2}\scal{P_U(v)}{v}\\
&= \scal{P_A(x)}{v}
+ \frac{\beta-\scal{P_A(x)}{v}}{\|P_U(v)\|^2}\|P_U(v)\|^2\\
&= \beta.
\end{align}
\end{subequations}
Hence $Q(x)\in B$ and we are done. 
\end{proof}

\begin{proposition} 
The assumption that $P_U(v)\neq 0$ implies 
\begin{align}
(A\cap B)- c &= U\cap \{v\}^\perp, \label{e:henry1}\\
\big((A\cap B)-c\big)^\perp &= U^\perp +\RR v \label{e:henry2}.
\end{align}
for every $c\in A\cap B$.
\end{proposition} 
\begin{proof}
By \cref{p:neq1}, $A\cap B\neq\varnothing$. 
Let $c\in A\cap B$. 
Then
\begin{equation}
(A\cap B)-c
= (A-c)\cap(B-c) = U \cap \{v\}^\perp,
\end{equation}
which is \cref{e:henry1}. 
Hence, using also 
\cite[Theorem~9.35 and Corollary~9.37]{Deutsch}, we see that 
\begin{equation}
\big((A\cap B)-c\big)^\perp 
= \big(U \cap \{v\}^\perp \big)^\perp
=\overline{U^\perp+\RR v} = U^\perp+\RR v.
\end{equation}
Therefore, \cref{e:henry2} is verified and we are done.
\end{proof}

\begin{theorem}
\label{t:5}
The assumption that $P_U(v)\neq 0$ implies $Q=P_{A\cap B}$.
\end{theorem}
\begin{proof}
Let $x\in X$.
By \cref{p:neq1}, 
\begin{equation}
\label{e:220616c}
Q(x)\in A\cap B.
\end{equation}
Using \cref{e:Q}, we have 
\begin{subequations}
\label{e:220616e}
\begin{align}
x - Q(x)
&= 
x - P_A(x) - \frac{\beta-\scal{P_A(x)}{v}}{\|P_U(v)\|^2}P_U(v)\\
&= 
\big(P_U(x)+P_{U^\perp}(x)\big) - \big(a_0+P_U(x)\big) - \frac{\beta-\scal{P_A(x)}{v}}{\|P_U(v)\|^2}(v-P_{U^\perp}v)\\
&= 
\bigg(P_{U^\perp}(x)-a_0 + \frac{\beta-\scal{P_A(x)}{v}}{\|P_U(v)\|^2}P_{U^\perp}(v) \bigg) 
- \frac{\beta-\scal{P_A(x)}{v}}{\|P_U(v)\|^2} v\\
&\in U^\perp + \RR v.
\end{align}
\end{subequations}
Now \cref{e:220616c}, \cref{e:henry2}, and \cref{e:220616e} yield
\begin{equation}
\label{e:220616d}
x- Q(x) \in \big((A \cap B)-Q(x)\big)^\perp.
\end{equation}
Combining \cref{e:220616c} and \cref{e:220616d},
we conclude that $P_{A\cap B}(x)=Q(x)$. 
\end{proof}

\begin{remark}
Note that \cref{t:5} and \cref{e:Q} imply \cref{t:main}\cref{t:main2}. 
\end{remark}

\section{A utility version of the main result and the case of two hyperplanes}

\label{sec:uti}

The analysis in the previous section was simplified because 
of our assumption that $\|v\|=1$.
It is worthwhile to record our results when 
we drop this normalization requirement.

\begin{theorem}[trichotomy]
\label{t:uti}
Let $A$ be as in \cref{e:A} 
and let $H$ be a hyperplane given by 
\begin{equation}
H := \menge{x\in X}{\scal{x}{c}=\gamma}, 
\end{equation}
where $c\in X\smallsetminus \{0\}$ and $\gamma\in\RR$.
Let $x\in X$. 
Then exactly one of the following holds:
\begin{enumerate}
\item \label{t:uti1}
$P_U(c) = 0$, $\scal{P_A(x)}{c}=\gamma$, $A\cap H\neq\varnothing$, and $P_{A\cap H}(x)=P_A(x)$. 
\item \label{t:uti2}
$P_U(c) = 0$, $\scal{P_A(x)}{c}\neq\gamma$, $A\cap H= \varnothing$, and $P_{\Fix(P_AP_H)}(x)=P_A(x)$. 
\item \label{t:uti3}
$P_U(c)\neq 0$, $A\cap H\neq\varnothing$, and 
$\displaystyle P_{A\cap H}(x)=  
P_A(x) + \frac{\gamma-\scal{P_A(x)}{c}}{\|P_U(c)\|^2}P_U(c)$. 
\end{enumerate}
\end{theorem}
\begin{proof}
Suppose that 
\begin{equation}
\label{e:220616f}
v = \frac{c}{\|c\|}
\;\;\text{and}\;\;
\beta = \frac{\gamma}{\|c\|}.
\end{equation}
Then $H = B$ (see \cref{e:B}). 
Note that 
\begin{equation}
P_U(v) = \frac{P_U(c)}{\|c\|},
\end{equation}
which shows that $P_U(v)=0$ $\Leftrightarrow$ $P_U(c)=0$. 

\cref{t:uti1}: 
This is clear from \cref{c:1}\cref{c:1.1}.
\cref{t:uti2}: 
Combine \cref{c:1}\cref{c:1.2} with \cref{e:altE}. 
\cref{t:uti3}: 
Combining \cref{t:5}, \cref{e:Q}, and \cref{e:220616f} yields
\begin{subequations}
\begin{align}
P_{A\cap B}(x)
&= 
P_A(x) + \frac{\beta-\scal{P_A(x)}{v}}{\|P_U(v)\|^2}P_U(v)\\
&= 
P_A(x) + \frac{{\gamma}/{\|c\|}-\scal{P_A(x)}{c/\|c\|}}{\|P_U(c/\|c\|)\|^2}P_U(c/\|c\|)\\
&= 
P_A(x) + \frac{\gamma-\scal{P_A(x)}{c}}{\|P_U(c)\|^2}P_U(c), 
\end{align}
\end{subequations}
as claimed.
\end{proof}

\begin{corollary}[two hyperplanes]
\label{c:twohyp}
Suppose that 
\begin{equation}
H_1 := \menge{x\in X}{\scal{x}{c_1}=\gamma_1}
\;\;\text{and}\;\;
H_2 := \menge{x\in X}{\scal{x}{c_2}=\gamma_2},
\end{equation}
where $c_1,c_2$ lie in $X\smallsetminus\{0\}$,
and $\gamma_1,\gamma_2$ belong to $\RR$. 
Let $x\in X$. 
Then the following hold:
\begin{enumerate}
\item \label{c:twohyp1}
If $\scal{c_1}{c_2}^2 = \|c_1\|^2\|c_2\|^2$
and 
$\|c_1\|^2(\scal{x}{c_2}-\gamma)=\scal{c_1}{c_2}
(\scal{x}{c_1}-\gamma_1)$, 
then $H_1=H_2$ and 
\begin{equation}
P_{H_1\cap H_2}(x) = P_{H_1}(x) = x - \frac{\scal{x}{c_1}-\gamma_1}{\|c_1\|^2}c_1. 
\end{equation}
\item \label{c:twohyp2}
If $\scal{c_1}{c_2}^2 = \|c_1\|^2\|c_2\|^2$ and 
$\|c_1\|^2(\scal{x}{c_2}-\gamma)\neq \scal{c_1}{c_2}
(\scal{x}{c_1}-\gamma_1)$, 
then $H_1$ and $H_2$ are parallel but distinct ($H_1\cap H_2 = \varnothing$), and 
\begin{equation}
P_{\Fix(P_{H_1}P_{H_2})}(x) = P_{H_1}(x) = x - \frac{\scal{x}{c_1}-\gamma_1}{\|c_1\|^2}c_1. 
\end{equation}
\item \label{c:twohyp3}
If $\scal{c_1}{c_2}^2 \neq \|c_1\|^2\|c_2\|^2$, 
then $H_1\cap H_2\neq \varnothing$ and 
\begin{subequations}
\begin{align}
P_{H_1\cap H_2}(x) &= x 
+ \frac{\|c_2\|^2\big(\gamma_1-\scal{x}{c_1}\big)+\scal{c_1}{c_2}\big(\scal{x}{c_2}-\gamma_2\big)}{\|c_1\|^2\|c_2\|^2-\scal{c_1}{c_2}^2} c_1 \\
&\qquad + 
\frac{\|c_1\|^2\big(\gamma_2-\scal{x}{c_2}\big)+\scal{c_1}{c_2}\big(\scal{x}{c_1}-\gamma_1\big)}{\|c_1\|^2\|c_2\|^2-\scal{c_1}{c_2}^2} c_2. 
\end{align}
\end{subequations}
\end{enumerate}
\end{corollary}
\begin{proof}
We apply \cref{t:uti} with 
$A = H_1$, $H= H_2$, $c=c_2$, and $\gamma =\gamma_2$. 
We have 
\begin{equation}
\label{e:two1}
P_A(x) = P_{H_1}(x)
= x- \frac{\scal{x}{c_1}-\gamma_1}{\|c_1\|^2}c_1,
\end{equation}
$U = \{c_1\}^\perp$, and so 
\begin{equation}
\label{e:two2}
P_U(c) = P_{\{c_1\}^\perp}(c_2) = c_2 - P_{\RR c_1}(c_2)
= c_2 - \frac{\scal{c_2}{c_1}}{\|c_1\|^2}c_1. 
\end{equation}
Therefore, \cref{e:two2} implies
\begin{subequations}
\label{e:two3}
\begin{align}
\|P_U(c)\|^2
&= \bigg\|c_2-\frac{\scal{c_1}{c_2}}{\|c_1\|^2}c_1 \bigg\|^2\\
&= \|c_2\|^2 - \frac{2\scal{c_1}{c_2}^2}{\|c_1\|^2}
+ \frac{\scal{c_1}{c_2}^2}{\|c_1\|^4}\|c_1\|^2\\
&= 
\|c_2\|^2 - \frac{\scal{c_1}{c_2}^2}{\|c_1\|^2}. 
\end{align}
\end{subequations}
Hence 
\begin{equation}
\label{e:two4}
P_U(c)=0
\;\;\Leftrightarrow\;\;
\|P_U(c)\|^2=0
\;\;\Leftrightarrow\;\;
\scal{c_1}{c_2}^2 = \|c_1\|^2\|c_2\|^2. 
\end{equation}
Next, 
\cref{e:two1} implies 
\begin{align}
\label{e:two5}
\scal{P_{A}(x)}{c}
&= 
\scal{P_{H_1}(x)}{c_2}
= 
\scal{x}{c_2} - \frac{\scal{x}{c_1}-\gamma_1}{\|c_1\|^2}\scal{c_1}{c_2}.
\end{align}
Thus
\begin{equation}
\scal{P_A(x)}{c}=\gamma
\;\;\Leftrightarrow\;\;
\|c_1\|^2\big(\scal{x}{c_2}-\gamma_2\big)=\scal{c_1}{c_2}
\big(\scal{x}{c_1}-\gamma_1\big). 
\end{equation}

\cref{c:twohyp1}:
The hypothesis in this case 
corresponds to $P_U(c)=0$ and $\scal{P_A(x)}{c}=\gamma$. 
By \cref{t:uti}\cref{t:uti1}, $A\cap H\neq\varnothing$ which means 
the gap vector $g=0$  by \cref{c:1}\cref{c:1.1}. 
By \cref{p:1}, $H_1 \subseteq H_2$.
Because both $H_1$ and $H_2$ are hyperplanes, we have $H_1=H_2=H_1\cap H_2$.

\cref{c:twohyp2}:
The hypothesis in this case 
corresponds to $P_U(c)=0$ and $\scal{P_A(x)}{c}\neq \gamma$. 
The conclusion now follows from 
\cref{t:uti}\cref{t:uti2}. 

\cref{c:twohyp3}:
Using \cref{e:two4}, \cref{t:uti}\cref{t:uti3}, \cref{e:two1}, \cref{e:two5}, 
\cref{e:two3}, and \cref{e:two2},
we have $H_1\cap H_2 = A\cap H\neq\varnothing$ and
\begin{subequations}
\begin{align}
P_{H_1\cap H_2}(x)
&= x- \frac{\scal{x}{c_1}-\gamma_1}{\|c_1\|^2}c_1
+ \frac{\gamma_2 - \Big(\scal{x}{c_2} - \frac{\scal{x}{c_1}-\gamma_1}{\|c_1\|^2}\scal{c_1}{c_2}\Big)}{\|c_2\|^2 - \frac{\scal{c_1}{c_2}^2}{\|c_1\|^2}}\Big(c_2 - \frac{\scal{c_2}{c_1}}{\|c_1\|^2}c_1  \Big)\\
&= x- \frac{\scal{x}{c_1}-\gamma_1}{\|c_1\|^2}c_1\\
&\qquad 
+ \frac{\|c_1\|^2\big(\gamma_2 - \scal{x}{c_2}\big) + 
\scal{c_1}{c_2}\big(\scal{x}{c_1}-\gamma_1\big)}{\|c_1\|^2\|c_2\|^2-\scal{c_1}{c_2}^2}
\Big(c_2 - \frac{\scal{c_2}{c_1}}{\|c_1\|^2}c_1  \Big)\\
&= x 
 + \frac{\|c_1\|^2\big(\gamma_2 - \scal{x}{c_2}\big) + 
\scal{c_1}{c_2}\big(\scal{x}{c_1}-\gamma_1\big)}{\|c_1\|^2\|c_2\|^2-\scal{c_1}{c_2}^2}
c_2\\
& \qquad - \Big(\frac{\scal{x}{c_1}-\gamma_1}{\|c_1\|^2} 
+ \frac{\|c_1\|^2\big(\gamma_2 - \scal{x}{c_2}\big) + 
\scal{c_1}{c_2}\big(\scal{x}{c_1}-\gamma_1\big)}{\|c_1\|^2\|c_2\|^2-\scal{c_1}{c_2}^2}\frac{\scal{c_2}{c_1}}{\|c_1\|^2}\Big)c_1\\
&= x 
 + \frac{\|c_1\|^2\big(\gamma_2 - \scal{x}{c_2}\big) + 
\scal{c_1}{c_2}\big(\scal{x}{c_1}-\gamma_1\big)}{\|c_1\|^2\|c_2\|^2-\scal{c_1}{c_2}^2}
c_2\\
& \qquad + 
\frac{\|c_2\|^2\big(\gamma_1-\scal{x}{c_1}\big)+\scal{c_1}{c_2}\big(\scal{x}{c_2}-\gamma_2\big)}{\|c_1\|^2\|c_2\|^2-\scal{c_1}{c_2}^2} c_1,
\end{align}
\end{subequations}
as announced. 
\end{proof}

\begin{remark}
Having just found a formula for the projection onto the intersection of two hyperplanes, it is in principle possible to present a formula
for the intersection of three (or more) hyperplanes;
however, the result would of course be significantly more
complicated than the formulas presented in \cref{c:twohyp}. 
\end{remark}

We conclude this section with the following 
limiting example which shows that 
there does not appear to exist a straightforward extension of 
the main result in \cref{t:main}.
Indeed, \cref{ex:cone} below 
verifies that \cref{t:main}
does not generalize when we replace $A$ by a
cone $K$.
Observe that in this case $a_0=0$,
hence $U$ is replaced by $K$ as well.
\begin{example}
\label{ex:cone}
Suppose that $X=\RR^2$,
and that 
$K=\RR^2_{+}\coloneqq\menge{(\xi_1,\xi_2)\in \RR^2}{\xi_1\ge 0, \xi_2\ge 0}$. Let $v_1=\tfrac{1}{\sqrt{2}}(-1,-1)$,
let $v_2=\tfrac{1}{\sqrt{2}}(1,1)$, let $\beta_1=0$ and let $\beta_2=\tfrac{1}{\sqrt{2}}$.
Set $(\forall i\in \{1,2\})$
$B_i\coloneqq\menge{x\in \RR^2}{\scal{x}{v_i}=\beta_i}$.
Then the  following hold:
\begin{enumerate}
\item 
Let $ x\in K\smallsetminus\{(0,0)\}$.
Then: 
\begin{enumerate}
\item
\label{ex:cone:i:a}
$K\cap B_1=\{(0,0)\}$.
\item
\label{ex:cone:i:b}
$P_K(v_1)=(0,0)$.
\item 
\label{ex:cone:i:c}
$  P_{K\cap B_1}\equiv (0,0)$.
\item
\label{ex:cone:i:d}
$(0,0) = P_{K\cap B_1}(x)\neq P_K(x)=x$.
\end{enumerate}
\item 
Set 
$ S\coloneqq\menge{(\xi,0)\in \RR^2}{\xi>\sqrt{2}}$
 and set
 $ \widetilde{Q} \coloneqq P_K + \ds\frac{\beta_2-\scal{P_K(\cdot)}{v_2}}{\|P_K(v_2)\|^2}P_K(v_2)$.
Then: 
\begin{enumerate}
\item
\label{ex:cone:ii:a}
$K\cap B_2=\conv\{(1,0),(0,1)\}$.
\item
\label{ex:cone:ii:b}
$P_K(v_2)=v_2\neq (0,0)$.
\item
\label{ex:cone:ii:c}
$\widetilde{Q}=P_K+\big(\tfrac{1}{\sqrt{2}}-\scal{P_K(\cdot)}{v_2}\big)v_2$.
\item
\label{ex:cone:ii:d}
$(\forall (\xi,0)\in S)$ $P_{K\cap B_2}(\xi,0))=(1,0)$.
\item
\label{ex:cone:ii:e}
$(\forall (\xi,0)\in S)$ 
$\widetilde{Q}((\xi,0))=\tfrac{1}{2}(\sqrt{2}+\xi,\sqrt{2}-\xi)\not\in K$.
\end{enumerate}
\end{enumerate}
\end{example}

\section{A numerical experiment}

\label{sec:numexp}

In this section we provide a numerical experiment to 
evaluate the performance of the formula developed in 
\cref{c:twohyp} when employed to 
find the projection onto the intersection 
of finitely many hyperplanes.

We randomly generate 100 matrices $M$ each of size 
$10\times 50$. For each matrix $M$, we randomly generate 
$\overline{x}\in \RR^{50}$ and set $b=M\overline{x}$. This guarantees that 
the underdetermined system of equations 
$Mx=b$ is consistent, i.e.,
it has a solution. For each random instance of the matrix $M$ we randomly generate $100$ starting points.
Because the $i^{\rm th}$ row in each of the randomly generated systems of equation $Mx=b$ defines a hyperplane,
namely 
$H_i\coloneqq \menge{x\in \RR^{50}}{\scal{m_i}{x}=b_i}$, 
we  set 
\begin{equation}
P\coloneqq P_{H_{10}}P_{H_{9}}\cdots P_{H_2}P_{H_1}
\end{equation}
and 
\begin{equation}
Q\coloneqq
P_{H_{10}\cap H_9}P_{H_{8}\cap H_7}\cdots P_{H_2\cap H_1}.
\end{equation}

For each of the randomly generated problems with data $(M,b)$,
and for a randomly generated starting point $x_0$,
let $x^*= P_C(x_0)$, where $C = M^{-1}(b)$.
We generate then 
two sequences via 
\begin{equation}
(\forall\nnn)\quad 
p_n\coloneqq P^n x_0\;\;\text{and}\;\;
q_n\coloneqq Q^n x_0.
\end{equation}
Both $(p_n)_\nnn$ and $(q_n)_\nnn$ are incarnations of the method of \emph{cyclic projections}, and thus both sequences converge to $x^*$, by e.g., \cite[Theorem~9.27]{Deutsch} or \cite[Chapter~3]{ER}. 
At each iteration index $n$, we measure the decibel (dB)
value of the proximity function which we choose to be the relative distance of the iterate 
to the solution $x^*$: 
\begin{equation}
20\log_{10}\frac{\norm{p_n-x^*}}{\norm{p_0-x^*}}
\;\;\text{and}\;\;
20\log_{10}\frac{\norm{q_n-x^*}}{\norm{q_0-x^*}}.
\end{equation}

\begin{figure}[H] 
\begin{center}
\includegraphics[scale=0.85]{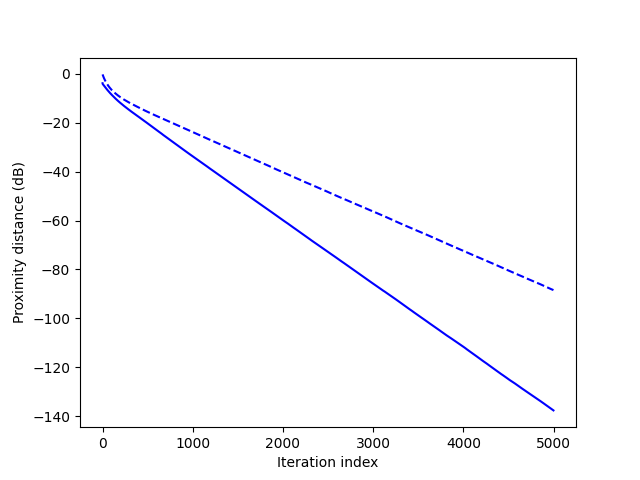}
\end{center}	
\caption{Plot of the decibel (dB)
value of the median of the proximity function
for the sequences $(p_n)_{\nnn}$ (the dashed curve)
and $(q_n)_{\nnn}$ (the solid curve).}
\label{figure}
\end{figure}

\Cref{figure} reports the progress of the proximity function of both sequences 
as a function of the iteration index where the median is calculated over all 100 instances of the matrix 
$M$ and then over $100$ randomly generated starting points, 
resulting in $10,000$ numerical scenarios. 
We observe a notable improvement in the speed of convergence using 
\cref{c:twohyp} --- this suggests that experimenting with this result may improve performance of projection algorithms involving hyperplanes. 

\section*{Declaration of competing interest}
\small
The authors declare that they have no competing interest.

\section*{Acknowledgements}
\small
The research of HHB and WMM was partially supported by Discovery Grants
of the Natural Sciences and Engineering Research Council of
Canada. The research of DM is partially supported by the 
University of Waterloo Mathematics Undergraduate Research Award.

\end{document}